\documentclass[preprint,10pt]{elsarticle}
\usepackage{amsfonts}

 \usepackage{bbm}
 \usepackage{amsmath}
 \usepackage{amssymb}
 \usepackage{mathrsfs}
 \usepackage[all]{xy}
 \usepackage{amsthm,amsmath,amsfonts,amssymb}
\usepackage{graphicx}

\biboptions{numbers,sort&compress}

\newtheorem{lem}{Lemma}[section]

\newtheorem{Acknw}{Acknowledgement}

\theoremstyle{definition}

\def\beq{\begin{equation}}
\def\eeq{\end{equation}}
\def\ben{\begin{equation*}}
\def\een{\end{equation*}}
\def \b{\begin}
\def\e{\end}

\begin{document}
\numberwithin{equation}{section}

\begin{frontmatter}
\title{Subsolution theorem and the Dirichlet problem for the quaternionic Monge-Amp\`{e}re equation\tnoteref{mytitlenote}}
\tnotetext[mytitlenote]{This work is supported by National Nature Science Foundation in China (No. 11571305).}

\author{Dongrui Wan\corref{mycorrespondingauthor}}
\cortext[mycorrespondingauthor]{Corresponding author}
\ead{wandongrui@szu.edu.cn,+86 755 2653 8953,
Room 415, Science and Technology Building, Shenzhen University}

\address{College of Mathematics and Statistics, Shenzhen University, Room 415 Science and Technology Building, Shenzhen, 518060, PR. China}

\begin{abstract}In this paper, the author studies quaternionic Monge-Amp\`{e}re equations and obtain the existence of the solutions to the Dirichlet problem for such equations in strictly pesudoconvex domains in quaternionic space. The stability and subsolution theorems are established for quaternionic Monge-Amp\`{e}re equations.
\end{abstract}


\end{frontmatter}
\section{Introduction}
In this paper, we are interested in solving, under possibly weak assumptions on the measure $d\mu$, the following Dirichlet problem for the quaternionic Monge-Amp\`{e}re equation in a given strictly pesudoconvex domain $\Omega$ in $n$-dimensional quaternionic space $\mathbb{H}^n$:
\beq\label{d4}\begin{cases}
 \,\,\,u \in PSH\cap L^\infty(\Omega)\\
 \,\,\,(\triangle u)^n= d\mu\\
 \,\,\,\lim_{\zeta\rightarrow q}u(\zeta)=\varphi(q)\,\,\,\,\,q\in\partial\Omega,~\varphi\in C(\partial\Omega),
 \end{cases}\eeq
where $(\triangle u)^n$ denotes the quaternionic Monge-Amp\`{e}re measure of $u$. We have shown in \cite{wan3} that quaternionic Monge-Amp\`{e}re operator $(\triangle u)^n$ is well defined as a positive measure for locally bounded quaternionic plurisubharmonic ($PSH$, for short) function $u$.

The quaternionic Monge-Amp\`{e}re equation, relating to the quaternionic version of Calabi-Yau conjucture, has attracted many analysts to study on it. Quaternionic analysis has important applications in the supersymmetric theory in physics.  It is interesting to study the quaternionic Monge-Amp\`ere operator over quaternionic manifolds, in particular to study the Dirichlet problem for the quaternionic Monge-Amp\`{e}re equation \cite{alesker7,alesker2,alesker6,alesker2017,alesker8,wang1,wang18,zhu}. 

The Dirichlet problem (\ref{d4}) in smooth case, for a smooth positive measure on the right hand side, i.e. $d\mu=gdV$ and $g>0$ smooth (denote by dV the Lebesgue measure), was first considered by Alesker \cite{alesker4} on a Euclidean ball $B$ in $\mathbb{H}^n$, then was recently solved by Zhu \cite{zhu} on general domains in $\mathbb{H}^n$. For continuous data, Alesker \cite{alesker4} showed that the unique continuous solution exists for $d\mu=gdV$, $0\leq g\in C(\overline{\Omega})$, and for strictly pesudoconvex domain $\Omega$ in $\mathbb{H}^n$. The existence of the solutions to the Dirichlet problem (\ref{d4}) for the quaternionic Monge-Amp\`{e}re equation when the right hand side $\mu$ is some more general positive Borel measure is still an open problem.

The existence theorems for the complex Monge-Amp\`{e}re equation have been undergoing intensive research in the past few decades. In mid-seventies Bedford and Taylor \cite{beddirichlet,bed} solved the following Dirichlet problem for complex Monge-Amp\`{e}re equation for $d\mu=gdV, g\in C(\overline{\Omega})$, in a strictly pesudoconvex domain in $\mathbb{C}^n$:
\beq\label{d5}\begin{cases}
 \,\,\,u \in PSH\cap L^\infty(\Omega)\\
 \,\,\,(dd^c u)^n= d\mu\\
 \,\,\,\lim_{\zeta\rightarrow z}u(\zeta)=\varphi(z)\,\,\,\,\,z\in\partial\Omega,~\varphi\in C(\partial\Omega).
 \end{cases}\eeq Cegrell \cite{cegrell84} generalized this result to the case of bounded $g$. Cegrell and Persson \cite{cp} also showed that continuous solutions exist if $d\mu=gdV, g\in L^2(\Omega,dV)$, but for $g\in L^1(\Omega,dV)$ this is not necessarily true \cite{cs}. Ko\l odziej \cite{kolodziej98,kolodziej99} proved that the above Dirichlet problem still admits a unique weak continuous solution when the right hand side $\mu$ is a measure satisfying some sufficient condition. Since this sufficient condition is not easy to varify, Ko\l odziej obtained the subsolution theorem, saying that the Dirichlet problem (\ref{d5}) in a strictly pesudoconvex domain is solvable if there is a subsolution (cf. \cite{kolodziej94,kolodziej95}). We refer to \cite{kolodziej05} for a nice survey on this.

Inspired by the idea of Cegrell \cite{cp} we use the connection between real and quaternionic Monge-Amp\`ere measure, which we established in \cite{wan6},  to obtain that there exists a unique solution to the Dirichlet problem (\ref{d4}) if $d\mu=gdV, g\in L^4(\Omega,dV)$. Let $|\cdot|_\Omega$ and $|\cdot|_{\partial\Omega}$ denote the sup-norm on $\Omega$ and $\partial\Omega$. 

\b{thm}\label{t1.2}Let $\Omega$ be a strictly pesudoconvex domain in $\mathbb{H}^n$. If $0\leq g\in L^4(\Omega)$, $\varphi\in C(\partial\Omega)$, then the Dirichlet problem (\ref{d4}) with $d\mu=gdV$ has a unique solution. The solution, denoted by $U_{\mathbb{Q}}(\varphi,g)$, is in $C(\overline{\Omega})$ and satisfies
 $$\inf_{\partial\Omega} \varphi-C\|g\|_{L^4(\Omega)}^{\frac{1}{n}}\leq U_{\mathbb{Q}}(\varphi,g)\leq\sup_{\partial\Omega}\varphi$$ for some constant $C$ depending only on $\Omega$. And
 $$|U_{\mathbb{Q}}(\varphi_1,g_1)-U_{\mathbb{Q}}(\varphi_2,g_2)|_{\Omega}\leq |\varphi_1-\varphi_2|_{\partial \Omega}+C\|g_1-g_2\|_{L^4(\Omega)}^{\frac{1}{n}}$$ for $\varphi_1,\varphi_2\in C(\partial \Omega)$ and $0\leq g_1,g_2\in L^4(\Omega)$.
\e{thm}
As a direct concequence of Theorem \ref{t1.2}, we obtain the $L^\infty-L^4-$stability of the quaternionic Monge-Amp\`ere equation. By H\"older's inequality we also get the $L^p-L^q-$stability of the quaternionic Monge-Amp\`ere equation for 
$$(p,q)\in ([1,n]\times[1,\infty])\cup ([1,\infty]\times[4,\infty]).$$

In this paper we show that the subsolution theorem for quaternionic Monge-Amp\`ere equation is still true by combining the above stability theorem and the well known method used by Ko\l odziej (cf. \cite{kolodziej94,kolodziej95,kolodziej05}).

\b{thm}\label{t1.5} Let $\Omega$ be a strictly pesudoconvex domain in $\mathbb{H}^n$. Let $\mu$ be a finite positive Borel measure in $\Omega$ and $\varphi\in C(\partial\Omega)$. If there exists a subsolution $v$, i.e.
\ben\label{d3}\begin{cases}
 \,\,\,v \in PSH\cap L^\infty(\Omega)\\
 \,\,\,(\triangle v)^n\geq d\mu\\
 \,\,\,\lim_{\zeta\rightarrow q}v(\zeta)=\varphi(q)\,\,\,\,\text{for~any}\,\,q\in\partial\Omega,
 \end{cases}\een then there exists a solution $u$ of the Dirichlet problem (\ref{d4}).
\e{thm}

Historically, the quaternionic Monge-Amp\`{e}re operator was  defined by Alesker \cite{alesker1,alesker4} as the Moore determinant of the quaternionic Hessian of $u$:
 \begin{equation}\label{det} det(u)=det\left[\frac{\partial^2u}{\partial \overline{q}_j\partial q_k}(q)\right].\end{equation} Compared to the complex pluripotential theory, the main difficulties in the quaternionic pluripotential theory come from the noncommutativity of elements of the quaternionic Hessian and the complexity of the noncommutative Moore determinant. Alesker \cite{alesker2,alesker2017} observed the relationship between the Baston operator and the quaternionic Monge-Amp\`{e}re operator. The Baston operator, which is the first operator of $0$-Cauchy-Fueter complex, is known explicitly \cite{baston,Wang,wang18}. Based on this observation, the author and Wang \cite{wan3} introduced the first order differential operators $d_0$ and $d_1$ and wrote the quaternionic Monge-Amp\`{e}re operator as the $n$-th exterior power of the Baston operator $\triangle=d_0d_1$. Then several results in the complex pluripotential theory were established for the quaternionic Monge-Amp\`{e}re operator \cite{wan4,wan6,wan7,wan8}. The theory of quaternionic closed positive currents established by the author and Wang \cite{wan3}, is an essential reason why we can still obtain the subsolution theorem for the quaternionic Monge-Amp\`{e}re equation as Ko\l odziej did for the complex case. 

\section{Preliminaries}
In this section, we are going to recall some basic definitions and conclusions of quaternionic Monge-Amp\`{e}re operator and quaternionic closed positive currents following \cite{wan3,wan4}.

A real valued function $f:\mathbb{H}^n\rightarrow \mathbb{R}$ is called \emph{quaternionic plurisubharmonic} ($PSH$, for short) if it is upper semi-continuous and its restriction to any right quaternionic line is subharmonic (in the usual sense). Any quaternionic $PSH$ function is subharmonic (cf. \cite{alesker1,alesker4,alesker2} for more information about $PSH$ functions). Denote by $PSH(\Omega)$ the class of quaternionic plurisubharmonic functions on $\Omega$.

\b{pro}\label{p1.1} Let $\Omega$ be an open subset of $\mathbb{H}^n$.\\
(1). The family $PSH(\Omega)$ is a convex cone, i.e. if $\alpha,\beta$ are non-negative numbers and $u,v\in PSH(\Omega)$, then $\alpha u+\beta v\in PSH(\Omega)$; and $\max\{u,v\}\in PSH(\Omega)$.\\
(2). If $\Omega$ is connected and $\{u_j\}\subset PSH(\Omega)$ is a decreasing sequence, then $u=\lim_{j\rightarrow\infty}u_j\in PSH(\Omega)$ or $u\equiv-\infty$.\\
(3). Let $\{u_\alpha\}_{\alpha\in A}\subset PSH(\Omega)$ be such that its upper envelope $u=\sup_{\alpha\in A}u_\alpha$ is locally bounded above. Then the upper semicontinuous regularization $u^*\in PSH(\Omega)$.\\
(4). Let $\omega$ be a non-empty proper open subset of $\Omega$, $u\in PSH(\Omega),v\in PSH(\omega)$, and $\limsup_{q\rightarrow\zeta}v(q)\leq u(\zeta)$ for each $\zeta\in \partial\omega\cap\Omega$, then
$$w:=\left\{
      \begin{array}{ll}
        \max\{u,v\}, & \text{in}~\omega \\
        u, &  \text{in}~\Omega\backslash\omega
      \end{array}
    \right. \quad\in PSH(\Omega).$$\\
(5). If $u\in PSH(\Omega)$, then the standard regularization $u_\varepsilon:=u*\rho_\varepsilon$ is also $PSH$ in $\Omega_\varepsilon:=\{q\in \Omega:dist(q,\partial\Omega)>\varepsilon\}$, moreover, $u_\varepsilon\searrow u$ as $\varepsilon\rightarrow0$.\\
(6). If $\gamma(t)$ is a convex increasing function of a parameter $t\in\mathbb{R}$ and $u\in PSH$, then $\gamma\circ u\in PSH$.
\e{pro}

We use the well known embedding of the quaternionic algebra $\mathbb{H}$ into End$(\mathbb{C}^{2})$ defined by
\begin{equation*}x_0+x_1\textbf{i}+x_2\textbf{j}+x_3\textbf{k}\mapsto\left(
                                                                                 \begin{array}{cc}
                                                                                   x_0+\textbf{i}x_1 & -x_2-\textbf{i}x_3 \\
                                                                                   x_2-\textbf{i}x_3 & x_0-\textbf{i}x_1 \\
                                                                                 \end{array}
                                                                               \right).
\end{equation*}
Actually we will use the conjugate embedding
\begin{equation*}\begin{aligned}\tau:\mathbb{H}^{n}\cong\mathbb{R}^{4n}&\hookrightarrow\mathbb{C}^{2n\times2},\\ (q_0,\ldots,q_{n-1})&\mapsto \textbf{z}=(z^{j\alpha})\in\mathbb{C}^{2n\times2},
\end{aligned}\end{equation*}
$j=0,1,\ldots,2n-1, ~\alpha=0 ,1 ,$ with
\begin{equation*}\label{2.2}\left(
                             \begin{array}{cc}
                               z^{(2l)0 } & z^{(2l)1 } \\
                               z^{(2l+1)0 } & z^{(2l+1)1 } \\
                             \end{array}
                           \right):=\left(
                                      \begin{array}{cc}
                                         x_{4l}-\textbf{i}x_{4l+1} & -x_{4l+2}+\textbf{i}x_{4l+3} \\
                                        x_{4l+2}+\textbf{i}x_{4l+3} & x_{4l}+\textbf{i}x_{4l+1} \\
                                      \end{array}
                                    \right),
\end{equation*}for $q_l=x_{4l}+\textbf{i}x_{4l+1}+\textbf{j}x_{4l+2}+\textbf{k}x_{4l+3}$, $l=0,1,\ldots,n-1$. Pulling back to the quaternionic
space $\mathbb{H}^n\cong\mathbb{R}^{4n}$ by the embedding above, we define on $\mathbb{R}^{4n}$ first-order differential operators $\nabla_{j\alpha}$ as following:
\begin{equation}\label{2.4}\left(
                             \begin{array}{cc}
                               \nabla_{(2l)0 } & \nabla_{(2l)1 } \\
                               \nabla_{(2l+1)0 } & \nabla_{(2l+1)1 } \\
                             \end{array}
                           \right):=\left(
                                      \begin{array}{cc}
                                         \partial_{x_{4l}}+\textbf{i}\partial_{x_{4l+1}} & -\partial_{x_{4l+2}}-\textbf{i}\partial_{x_{4l+3}} \\
                                        \partial_{x_{4l+2}}-\textbf{i}\partial_{x_{4l+3}} & \partial_{x_{4l}}-\textbf{i}\partial_{x_{4l+1}} \\
                                      \end{array}
                                    \right).
\end{equation}
$z^{k\beta}$'s can be viewed as independent variables and $\nabla_{j\alpha}$'s are derivatives with respect to these variables. The operators $\nabla_{j\alpha}$'s play very important roles in the investigating of regular functions in several quaternionic variables \cite{Wang,kang}.
\par

Let $\wedge^{2k}\mathbb{C}^{2n}$ be the complex exterior algebra generated by $\mathbb{C}^{2n}$, $0\leq k\leq n$. Fix a basis
$\{\omega^0,\omega^1,\ldots$, $\omega^{2n-1}\}$ of $\mathbb{C}^{2n}$. Let $\Omega$ be a domain in $\mathbb{R}^{4n}$. Define $d_0,d_1:C_0^\infty(\Omega,\wedge^{p}\mathbb{C}^{2n})\rightarrow C_0^\infty(\Omega,\wedge^{p+1}\mathbb{C}^{2n})$ by \begin{equation}\begin{aligned}\label{2.228}&d_0F=\sum_{k,I}\nabla_{k0 }f_{I}~\omega^k\wedge\omega^I,\\
&d_1F=\sum_{k,I}\nabla_{k1 }f_{I}~\omega^k\wedge\omega^I,\\
&\triangle
F=d_0d_1F,
\end{aligned}\end{equation}for $F=\sum_{I}f_{I}\omega^I\in C_0^\infty(\Omega,\wedge^{p}\mathbb{C}^{2n})$,  where the multi-index
$I=(i_1,\ldots,i_{p})$ and
$\omega^I:=\omega^{i_1}\wedge\ldots\wedge\omega^{i_{p}}$. Although $d_0,d_1$ are not exterior differential, their behavior is similar to the exterior differential.
\begin{lem}\label{p1.1}$($1$)$ $d_0d_1=-d_1d_0$. $($2$)$ $d_0^2=d_1^2=0$, thus $d_0\triangle=d_1\triangle=0$.\\
$($3$)$ For $F\in C_0^\infty(\Omega,\wedge^{p}\mathbb{C}^{2n})$, $G\in C_0^\infty(\Omega,\wedge^{q}\mathbb{C}^{2n})$, we have\begin{equation*}d_\alpha(F\wedge G)=d_\alpha F\wedge G+(-1)^{p}F\wedge d_\alpha G,\qquad \alpha=0,1.\end{equation*}
\end{lem}

We say $F$
is \emph{closed} if $d_0F=d_1F=0.$
For $u_1,\ldots,
u_n\in C^2$, $\triangle
u_1\wedge\ldots\wedge\triangle u_k$ is closed, $k=1,\ldots,n $. Moreover, it follows easily from Lemma \ref{p1.1} that $\triangle u_1\wedge\ldots\wedge\triangle u_n$ satisfies the following remarkable identities:\begin{equation}\label{2.37}\begin{aligned}\triangle u_1\wedge \triangle
u_2\wedge\ldots\wedge\triangle u_n&=d_0(d_1u_1\wedge \triangle
u_2\wedge\ldots\wedge\triangle u_n)\\&=d_0d_1(u_1\triangle
u_2\wedge\ldots\wedge\triangle u_n)\\&=\triangle (u_1
\triangle u_2\wedge\ldots\wedge\triangle u_n).
\end{aligned}\end{equation}

To write down the explicit expression, we
define for a function $u\in C^2$,\begin{equation*}\label{2.10}\triangle_{ij}u:=\frac{1}{2}(\nabla_{i0 }\nabla_{j1 }u-\nabla_{i1 }\nabla_{j0 }u),~~~\triangle u=\sum_{i,j=0}^{2n-1}\triangle_{ij}u~\omega^i\wedge
\omega^j.\end{equation*}
\begin{equation}\label{2.11}\begin{aligned}\triangle
u_1\wedge\ldots\wedge\triangle
u_n&=\sum_{i_1,j_1,\ldots}\triangle_{i_1j_1}u_1\ldots\triangle_{i_nj_n}u_n~\omega^{i_1}\wedge
\omega^{j_1}\wedge\ldots\wedge \omega^{i_n}\wedge
\omega^{j_n}\\&=\sum_{i_1,j_1,\ldots}\delta^{i_1j_1\ldots
i_nj_n}_{01\ldots(2n-1)}\triangle_{i_1j_1}u_1\ldots\triangle_{i_nj_n}u_n~\Omega_{2n},
\end{aligned}\end{equation}for $u_1,\ldots,u_n\in
C^2$, where $\Omega_{2n}$ is defined as\begin{equation}\label{2.21}\Omega_{2n}:=\omega^0\wedge
\omega^1\wedge\ldots\wedge\omega^{2n-2}\wedge
\omega^{2n-1},\end{equation}and $\delta^{i_1j_1\ldots
i_nj_n}_{01\ldots(2n-1)}:=$ the sign of the permutation from $(i_1,j_1,\ldots
i_n,j_n)$ to $(0,1,\ldots,2n-1),$ if $\{i_1,j_1,\ldots,
i_n,j_n\}$ $=\{0,1,\ldots,2n-1\}$; otherwise, $\delta^{i_1j_1\ldots
i_nj_n}_{01\ldots(2n-1)}=0$.

Note that $\triangle
u_1\wedge\ldots\wedge\triangle u_n$ is symmetric with respect to the
permutation of $u_1,\ldots,u_n$. In particular, when
$u_1=\ldots=u_n=u$, $\triangle u_1\wedge\ldots\wedge\triangle u_n$
coincides with $(\triangle u)^n:=\wedge^n\triangle u=\triangle_nu~\Omega_{2n}$. Denote by $\triangle_n(u_1,\ldots, u_n)$ the
coefficient of the form $\triangle u_1\wedge\ldots\wedge\triangle
u_n$, i.e., $\triangle u_1\wedge\ldots\wedge\triangle
u_n=\triangle_n(u_1,\ldots, u_n)~\Omega_{2n}.$ Then
$\triangle_n(u_1,\ldots,
u_n)$ coincides with the mixed
Monge-Amp\`{e}re operator $\text{det}(u_1,\ldots,u_n)$ while $\triangle_nu$
coincides with the quaternionic Monge-Amp\`{e}re operator
$\text{det}(u)$. This result was proved by
Alesker (Proposition 7.1 in \cite{alesker2}). And in Appendix A in \cite{wan3} we gave  an elementary and simpler proof of the identity:\beq\label{identity}\triangle u_1\wedge\ldots\wedge\triangle u_n=n!~\text{det}~(u_1,\ldots,u_n)\Omega_{2n}\eeq for $C^2$ function $u$.

By introducing the quaternionic version of differential forms, the author and Wang defined in \cite{wan3} the notions of closed positive forms and closed positive currents in the quaternionic case and our definition of closedness matches positivity well. Although a $2n$-form is not an authentic differential form and we cannot integrate it, we can define $\int_\Omega F:=\int_\Omega f dV,$ if we write $F=f~\Omega_{2n}\in L^1(\Omega,\wedge^{2n}\mathbb{C}^{2n})$,
where $dV$ is the Lebesgue measure and $\Omega_{2n}$ is given by (\ref{2.21}). In particular, if $F$ is positive $2n$-form, then $\int_\Omega F
\geq0$. For a $2n$-current $F=\mu~\Omega_{2n}$ with coefficient to be measure $\mu$, define
\begin{equation}\label{integral}\int_\Omega F:=\int_\Omega \mu.\end{equation}

\begin{lem}\label{p3.9}$($Stokes-type formula, Lemma 3.2 in \cite{wan3} $)$ Assume that $T=\sum_iT_i\omega^{\widehat{i}}$ is a smooth $(2n-1)$-form in $\Omega$, where $\omega^{\widehat{i}}=\omega^0\wedge\ldots\wedge\omega^{i-1}\wedge\omega^{i+1}\wedge\ldots\wedge\omega^{2n-1}$. Then for smooth function $h$, we have \begin{equation*}\int_\Omega hd_\alpha T=-\int_\Omega d_\alpha h\wedge T+\sum_{i=0}^{2n-1}(-1)^{i-1}\int_{\partial\Omega}hT_i~n_{i\alpha}~ dS,\end{equation*}where $ n_{i\alpha}$, $i=0,1,\ldots,2n-1$, $\alpha=0,1$, is defined by the matrix:\begin{equation}\label{sss4.3}\begin{aligned}\left(
                                                                                                    \begin{array}{cc}
                                                                                                      n_{(2l)0 }& n_{(2l)1 }\\
                                                                                                      n_{(2l+1)0 }&n_{(2l+1)1 }\\
                                                                                                    \end{array}
                                                                                                  \right):=&\left(
                                                                                                                                  \begin{array}{cc}
                                                                                                                                    n_{4l}+\emph{\textbf{i}}~n_{4l+1}&-n_{4l+2}-\emph{\textbf{i}}~n_{4l+3}\\
                                                                                                                                    n_{4l+2}-\emph{\textbf{i}}~n_{4l+3}&n_{4l}-\emph{\textbf{i}}~n_{4l+1}\\
                                                                                                                                  \end{array}
                                                                                                                                \right),
\end{aligned}\end{equation}$l=0,1,\ldots,n-1$. Here \emph{\textbf{n}}$=(n_0,n_1,\ldots,n_{4n-1})$ is the unit outer normal vector to $\partial\Omega$ and $dS$ denotes the surface measure of $\partial\Omega$. In particular, if $h=0$ on $\partial\Omega$, we have \begin{equation}\label{stokes}\int_\Omega hd_\alpha T=-\int_\Omega d_\alpha h\wedge T,\qquad \alpha=0,1.\end{equation}
\end{lem}

Bedford-Taylor theory \cite{bed,bed1990} in complex analysis can be generalized to the quaternionic case. Let $u$ be a locally bounded $PSH$ function and let $T$ be a closed
positive $2k$-current. Then $\triangle u\wedge T$ defined by \begin{equation}\label{3.110}\triangle u\wedge
T:=\triangle(uT),\end{equation}i.e., $(\triangle u\wedge
T)(\eta):=uT(\triangle\eta)$, for any test form $\eta$, is also a closed positive
current. Inductively,
\begin{equation}\label{3.111}\triangle
u_1\wedge\ldots\wedge\triangle u_p\wedge T:=\triangle(u_1\triangle
u_2\wedge\ldots\wedge\triangle u_p\wedge T)\end{equation} is a closed
positive current, when $u_1,\ldots,u_p\in PSH\cap
L_{loc}^\infty(\Omega)$.
In particular, for $u_1,\ldots,u_n\in PSH\cap
L_{loc}^{\infty} (\Omega)$, $\triangle u_1\wedge\ldots\wedge\triangle u_n=\mu\Omega_{2n}$ for a well defined positive Radon
measure $\mu$.

For any strongly positive test $(2n -2p)$-form $\psi$ on $\Omega$, (\ref{3.111}) can be rewritten as
\begin{equation}\label{3.26}\int_\Omega\triangle
u_1\wedge\ldots\wedge\triangle u_p \wedge\psi=\int_\Omega u_1\triangle
u_2\wedge\ldots\wedge\triangle u_p \wedge\triangle\psi, \end{equation}where $u_1,\ldots,u_p\in PSH\cap
L_{loc}^\infty(\Omega)$. Since positive currents have measure coefficients, (\ref{3.26}) also holds for strongly positive $\psi\in\mathcal{D}_0^{2n-2p}(\Omega)$ vanishing on the boundary.

The following different types of weak convergence results are powerful tools in developing pluripotential theory for the quaternionic Monge-Amp\`{e}re operator. We will use these results frequently in the following.
\begin{lem}\label{l3.2}(1) (Theorem 3.1 in \cite{wan3}) \\Let
$v^1,\ldots,v^k\in PSH\cap L_{loc}^{\infty}(\Omega)$. Let
$\{v_j^1\}_{j\in\mathbb{N}},\ldots,\{v_j^k\}_{j\in\mathbb{N}}$ be
decreasing sequences of PSH functions in $\Omega$ such
that $\lim_{j\rightarrow\infty}v_j^t=v^t$   pointwisely in $\Omega$
for each $t $. Then the currents $\triangle
v_j^1\wedge\ldots\wedge\triangle v_j^k $ converge weakly to
$\triangle v^1\wedge\ldots\wedge\triangle v^k $ as
$j\rightarrow\infty$.\\
(2) (Proposition 3.2 in \cite{wan4})\\ Let $v^0,\ldots,v^k\in PSH\cap L_{loc}^\infty(\Omega)$. Let
$\{v_j^0\}_{j\in\mathbb{N}},\ldots,\{v_j^k\}_{j\in\mathbb{N}}$ be
decreasing sequences of PSH functions in $\Omega$ such
that $\lim_{j\rightarrow\infty}v_j^t=v^t$ pointwisely in $\Omega$
for $t=0,\ldots,k$. Then the currents $v_j^0\triangle
v_j^1\wedge\ldots\wedge\triangle v_j^k $ converge weakly to
$v^0\triangle v^1\wedge\ldots\wedge\triangle v^k $ as
$j\rightarrow\infty$.\\
(3) (Theorem 2.1.11 in \cite{alesker1}) Let $\{v_j\}$ be a sequence of continuous $PSH$ functions in $\Omega$. Assume that this sequence converges uniformly on compact subsets to a function $v$. Then $v$ is continuous $PSH$ function. Moreover the measures $(\triangle v_j)^n$ converge weakly to $(\triangle v)^n$ as
$j\rightarrow\infty$. \end{lem}

\section{Solvability and stability of quaternionic Monge-Amp\`ere equation}
In this section, we are to prove Theorem \ref{t1.2} by combining the well known results for the Dirichlet problem of real Monge-Amp\`ere equation and the connection between real and quaternionic Monge-Amp\`ere operator. We refer to \cite{caffarelli1,taylor} for more detailed historical discussions for the real Monge-Amp\`ere equation. Here we only mention the following two basic results.

Let $|\cdot|_\Omega$ and $|\cdot|_{\partial\Omega}$ denote the sup-norm on $\Omega$ and $\partial\Omega$. Denote by $det_{\mathbb{R}}u$ the real Monge-Amp\`ere measure of $u$ in the usual sense and denote by $det(u)$ the Moore determinant of quaternionic Hessian of $u$ given by (\ref{det}).

\b{lem}\label{l1.1}(Theorem 4.1 and Lemma 3.5 in \cite{taylor}) Let $\Omega$ be bounded strictly convex in $\mathbb{R}^m$. Let $\varphi \in C(\partial \Omega)$ and $0\leq g\in L^1(\Omega)$. Then the following Dirichlet problem has a unique solution:
\beq\label{d1}\begin{cases}
 \,\,\,u \text{~is~convex~in~} \Omega\\
 \,\,\,det_{\mathbb{R}}u=gdV\\
 \,\,\,\lim_{x\rightarrow \xi}u(x)=\varphi (\xi)\,\,\,\,\text{for~all}\,\,\,\,\xi\in\partial\Omega.\\
 \end{cases}\eeq
Furthermore, the solution, denoted by $U_{\mathbb{R}}(\varphi,g)$, satisfies
 $$\inf_{\partial\Omega} \varphi-C\|g\|_{L^1(\Omega)}^{\frac{1}{m}}\leq U_{\mathbb{R}}(\varphi,g)\leq\sup_{\partial\Omega}\varphi$$ for some constant $C$ depending only on $\Omega$. And
 $$|U_{\mathbb{R}}(\varphi_1,g_1)-U_{\mathbb{R}}(\varphi_2,g_2)|_{\Omega}\leq |\varphi_1-\varphi_2|_{\partial \Omega}+C\|g_1-g_2\|_{L^1(\Omega)}^{\frac{1}{m}}$$ for $\varphi_1,\varphi_2\in C(\partial \Omega)$ and $0\leq g_1,g_2\in L^1(\Omega)$.
\e{lem}

\b{lem}\label{l1.2}(Theorem 1.1 in \cite{caffarelli1}) Let $\Omega$ be bounded strictly convex in $\mathbb{R}^m$, $\partial\Omega\in C^\infty$. If $\varphi\in C^\infty(\partial\Omega)$, $0<g\in C^\infty(\overline{\Omega})$, then the solution $U_{\mathbb{R}}(\varphi,g)$ of (\ref{d1}) exists and is in $C^\infty(\overline{\Omega})$.
\e{lem}
Alesker \cite{alesker4} and Zhu \cite{zhu} studied the Dirichlet problem for the quaternionic Monge-Amp\`ere equation in terms of the original definition of the quaternionic Monge-Amp\`ere operator as a Moore determinant $det(u)$. As we introduced in previous sections, we \cite{wan3} showed that the original definition $det(u)$ coincides with our new definition $(\triangle u)^n$. Since we need to use the relationship between real and quaternionic Monge-Amp\`ere operator (see Lemma \ref{l1.6} below), in this section we use the original definition of quaternionic Monge-Amp\`ere operator as a Moore determinant $det(u)$. 

To prove Theorem \ref{t1.2}, it is equivalent to prove the  same conclusion for the solution $U_{\mathbb{Q}}(\varphi,g)$ of the Dirichlet problem:
\beq\label{d6}\begin{cases}
 \,\,\,u \in PSH\cap L^\infty(\Omega)\\
 \,\,\,det(u)= gdV\\
 \,\,\,\lim_{\zeta\rightarrow q}u(\zeta)=\varphi(q)\,\,\,\,\,q\in\partial\Omega,~\varphi\in C(\partial\Omega).
 \end{cases}\eeq

Recall that an open bounded domain $\Omega\subset \mathbb{H}^n$ with a smooth boundary $\partial\Omega$ is called \emph{strictly pseudoconvex} if for every point $q_0\in\partial \Omega$ there exists a neighborhood $\mathcal{O}$ and a smooth strictly psh function $h$ on $\mathcal{O}$ such that $\Omega\cap \mathcal{O}=\{h<0\}, h(q_0)=0,$ and $\nabla h(q_0)\neq 0$.
\b{lem}\label{l1.3}(Corollary 1.3 in \cite{zhu}) Let $\Omega$ be a quaternionic strictly pseudoconvex bounded domain in $\mathbb{H}^n$. If $0<g\in C^\infty(\overline{\Omega})$, $\varphi\in C^\infty(\partial\Omega)$, then the Dirichlet problem (\ref{d6}) has a unique solution. And the solution is in $C^\infty(\overline{\Omega})$.
\e{lem}

\b{lem}\label{l1.4}(Theorem 1.3 in \cite{alesker4}) Let $\Omega$ be a quaternionic strictly pseudoconvex bounded domain in $\mathbb{H}^n$. If $0\leq g\in C(\overline{\Omega})$, $\varphi\in C(\partial\Omega)$, then the Dirichlet problem (\ref{d6}) has a unique solution. And the solution is in $C(\overline{\Omega})$.
\e{lem}

Here we rewrite the following comparison principle by using $det(u)$. We \cite{wan4} established these results in terms of $(\triangle u)^n$ by using the theory of quaternionic closed positive currents. We will use the comparison principle frequently in the following.
\begin{lem}\label{comparison}(Comparison principle, Theorem 1.2 in \cite{wan4})\\ (1) Let $u,v\in PSH\cap L_{loc}^\infty(\Omega)$. If for any $\zeta\in\partial\Omega$,
$\liminf_{\zeta\leftarrow q\in\Omega}(u(q)-v(q))\geq0,$
then
\begin{equation*}\label{5.1}\int_{\{u<v\}}det(v)\leq\int_{\{u<v\}}det(u).
\end{equation*}(2) Under the assumptions of (1), the inequality $det(u)\leq det(v)$ implies $v\leq u$.\\
(3) If for any $\zeta\in\partial\Omega$, $\liminf_{\zeta\leftarrow q\in\Omega}u(q)=\liminf_{\zeta\leftarrow q\in\Omega}v(q)=0$ and $u\leq v$ in $\Omega$, then 
$$\int_{\Omega}det(v)\leq\int_{\Omega}det(u).$$
\end{lem}

\b{cor}\label{c1.1} $$U_{\mathbb{Q}}(\varphi_1,g_1)+U_{\mathbb{Q}}(\varphi_2,g_2)\leq U_{\mathbb{Q}}(\varphi_1+\varphi_2,g_1+g_2)$$
$$|U_{\mathbb{Q}}(\varphi_1,g_1)-U_{\mathbb{Q}}(\varphi_2,g_2)|\leq -U_{\mathbb{Q}}(-|\varphi_1-\varphi_2|,|g_1-g_2|).$$
\e{cor}
\proof By the comparison principle and superadditivity.\endproof

\b{lem}\label{l1.6}(Proposition 4.1 in \cite{wan6}) For a function $u\in PSH\cap C^2$  , we have the inequality: $$\left(det(\frac{\partial^2u}{\partial\overline{q}_j\partial q_k})\right)^{\frac{1}{n}}\geq 4\left(det_\mathbb{R}(\frac{\partial^2u}{\partial x_s\partial x_t})\right)^{\frac{1}{4n}}.$$
\e{lem}
\b{pro}\label{t1.1}Let $\Omega$ be a strictly convex bounded domain in $\mathbb{H}^n$. If $0\leq g\in C(\overline{\Omega})$, $\varphi\in C(\partial\Omega)$, then $$U_{\mathbb{R}}(\varphi,C_0g^4)\leq U_{\mathbb{Q}}(\varphi,g), ~~C_0=4^{-4n}.$$
\e{pro}
\b{proof} By identifying $\mathbb{H}^n$ with $\mathbb{R}^{4n}$ and using Lemma \ref{l1.1}, $U_{\mathbb{R}}(\varphi,C_0g^4)$ exists and $u:=U_{\mathbb{R}}(\varphi,C_0g^4)$ is convex. By Lemma \ref{l1.4}, $U_{\mathbb{Q}}(\varphi,g)$ exists. By comparison principle, it is sufficient to prove $det (u)\geq g$, i.e., \beq\label{0.1}\int\psi det ( u)\geq \int \psi g,\eeq for any $ 0\leq \psi\in C_0(\Omega)$.

 Without loss of generality, we can assume that $\partial \Omega\in C^\infty$. Take a sequence $\varphi_j\in C^\infty(\partial \Omega)$ converging uniformly to $\varphi$, and a sequence $0<g_j\in C^\infty(\overline{\Omega})$ converging uniformly to $g$. Let $u_j=U_{\mathbb{R}}(\varphi_j,C_0g_j^4)$. Then $u_j$ converges uniformly to $u$ as $j\rightarrow \infty$ by Lemma \ref{l1.1}. By weak convergence result (Lemma \ref{l3.2}), $$\int \psi det  (u_j)\longrightarrow \int \psi det  (u) ~~\text{and}~~\int \psi g_j\longrightarrow \int \psi g.$$
Therefore, it is sufficient to prove (\ref{0.1}) for $0<g\in C^\infty(\overline{\Omega})$, $\varphi\in C^\infty(\partial\Omega)$, and $\partial \Omega\in C^\infty$. Note that in this case $u\in C^\infty(\overline{\Omega})$ by Lemma \ref{l1.2}. It suffices to prove that $$C_0(det  u)^4\geq C_0g^4=det_{\mathbb{R}} u$$ for $u\in C^\infty(\overline{\Omega})$. Then the conclusion follows from Lemma \ref{l1.6}. 
\e{proof}

\emph{Proof of Theorem \ref{t1.2}} First assume that $\Omega$ is strictly convex. It is sufficient to prove the inequalities for $g\in C(\overline{\Omega})$. For $0\leq g\in L^4(\Omega)$, we can take a sequence $0\leq g_j\in C(\overline{\Omega})$ converging to $g$ in $L^4(\Omega)$. Then $U_{\mathbb{Q}}(\varphi,g_j)\in C(\overline{\Omega})$ by Lemma \ref{l1.4}. By the second inequality, $$|U_{\mathbb{Q}}(\varphi,g_j)-U_{\mathbb{Q}}(\varphi,g_k)|_{\Omega}\leq C\|g_j-g_k\|_{L^4(\Omega)}^{\frac{1}{n}}.$$ Thus $u_j:=U_{\mathbb{Q}}(\varphi,g_j)$ converges uniformly to a continuous function $u$. And $u=U_{\mathbb{Q}}(\varphi,g)$ by the Lemma \ref{l3.2} and comparison principle (Lemma \ref{comparison}).

Apply Lemma \ref{l1.1} with $m=4n$ and Proposition \ref{t1.1} to get $$\inf_{\partial\Omega} \varphi-C\|C_0g^4\|_{L^1(\Omega)}^{\frac{1}{4n}}\leq U_{\mathbb{R}}(\varphi,C_0g^4)\leq U_{\mathbb{Q}}(\varphi,g)\leq\sup_{\partial\Omega}\varphi,~~~~i.e.,$$ $$\inf_{\partial\Omega} \varphi-C'\|g\|_{L^4(\Omega)}^{\frac{1}{n}}\leq U_{\mathbb{Q}}(\varphi,g)\leq\sup_{\partial\Omega}\varphi.$$ The second inequality follows from the first inequality and Corollary \ref{c1.1}.

For strictly pesudoconvex $\Omega$, we can take a bounded strictly convex domain $\Omega'$ containing $\Omega$ and extend $|g_1-g_2|$ by zero to an $L^4$-function on $\Omega'$. Then the theorem follows. \qed
\vskip 3mm

It follows from Theorem \ref{t1.2} that for $0\leq g\in L^4(\Omega)$, $$ \|U_{\mathbb{Q}}(0,g)\|_{L^\infty(\Omega)}\leq C\|g\|_{L^4(\Omega)}^{\frac{1}{n}}.$$ We obtain the $L^\infty-L^4-$stability of the quaternionic Monge-Amp\`ere equation. Then we get that the quaternionic Monge-Amp\`ere equation is $L^p-L^q-$stable for $p\in [1,\infty]$ and $q\geq4$ by H\"older's inequality.

\b{rem}Here we are inspired by the article \cite{cp} writen Cegrell and Persson. They proved the $L^\infty-L^2-$stability of the complex Monge-Amp\`ere equation by using the connection between real and complex Monge-Amp\`ere operators (due to the idea of Cheng and Yau mentioned in \cite{bed1990}). Blocki \cite{blocki} obtained the $L^n-L^1-$stability after showing an estimate for the complex Monge-Amp\`ere operator. For the quaternionic Monge-Amp\`ere equation, Blocki's estimate was generalized to the quaternionic case (cf. Lemma 4.6 in \cite{wan-cegrell}), thus one can also get the $L^n-L^1-$stability for quaternionic Monge-Amp\`ere equation by following Blocki's method. By H\"older's inequality we obtain the $L^p-L^q-$stability for quaternionic Monge-Amp\`ere equation for 
$$(p,q)\in ([1,n]\times[1,\infty])\cup ([1,\infty]\times[4,\infty]).$$ We do not know the stability for other pairs $(p,q)$.
\e{rem}

As an application of Theorem \ref{t1.2}, we generalize an inequality for the mixed quaternionic Monge-Amp\`ere measures. We prove this inequality in the smooth case in Appendix and generalize it to the nonsmooth functions in the following proposition. As for the complex Monge-Amp\`ere measure, the nonsmooth version of this inequality has nontrivial applications \cite{beddirichlet,demailly1992,kolodziej2003,dinew}. 

\b{lem}For $u_1,u_2,\ldots,u_n\in PSH\cap C^2(\Omega)$, \beq\label{5.4} det(u_1,\ldots,u_n)\geq det(u_1)^{\frac{1}{n}}\cdots det(u_n)^{\frac{1}{n}}.\eeq 
\e{lem}
\proof See Appendix for the proof.\endproof

\b{pro} \label{t1.4} Let $0\leq g,f \in L^1(B)$, $u,v\in PSH\cap C(\overline{B})$ such that $$det(u) \geq fdV,~~ det(v)\geq gdV.$$ Then $$det(\underbrace{u,\ldots,u}_k,\underbrace{v,\ldots, v}_{n-k})\geq f^{\frac{k}{n}}g^{\frac{n-k}{n}}dV.$$
\e{pro}
\proof 
For smooth $u,v$, we already have $$det(\underbrace{u,\ldots,u}_k,\underbrace{v,\ldots, v}_{n-k})\geq det(u)^{\frac{k}{n}}\cdot det(v)^{\frac{n-k}{n}}\geq f^{\frac{k}{n}}g^{\frac{n-k}{n}}dV.$$
Now for $g,f \in L^4(B)$, take sequences $\{g_j\},\{f_j\}$ such that  \ben\b{aligned}&0<g_j\in C^\infty\rightarrow g~~\text{in}~~ L^4(B), \\&0<f_j\in C^\infty\rightarrow f~~\text{in}~~ L^4(B).\e{aligned}\een Take $\varphi_j,\psi_j$ smooth on $\partial B$ and $\varphi_j\rightarrow u,~\psi_j\rightarrow v$ uniformly on $\partial B$. By Lemma \ref{l1.3}, there exist $u_j,v_j\in PSH\cap C^\infty(\overline{B})$ such that 
\ben\begin{cases}
 \,\,\,det(u_j)=f_jdV\,\,\,\,\text{in}\,\,B\\
 \,\,\,u_j=\varphi_j\,\,\,\,\text{on}\,\,\partial B
 \end{cases}~~~~~~~\begin{cases}
 \,\,\,det(v_j)=g_jdV\,\,\,\,\text{in}\,\,B\\
 \,\,\,v_j=\psi_j\,\,\,\,\text{on}\,\,\partial B.
 \end{cases}\een
By Theorem \ref{t1.2}, $u_j,v_j$ converge uniformly to $u,v$ respectively. Then it follows from Lemma \ref{l3.2} that
\beq\label{1.4}\begin{aligned} det(\underbrace{u,\ldots,u}_k,\underbrace{v,\ldots, v}_{n-k})=& \lim_{j\rightarrow \infty}det(\underbrace{u_j,\ldots,u_j}_k,\underbrace{v_j,\ldots, v_j}_{n-k})\\ \geq& \lim_{j\rightarrow \infty}(f_j)^{\frac{k}{n}}(g_j)^{\frac{n-k}{n}}dV \geq f^{\frac{k}{n}}g^{\frac{n-k}{n}}dV.\end{aligned}\eeq

Now for $0\leq g,f \in L^1(B)$, take $ g_j\in L^4(B)\nearrow g$ and $f_j\in L^4(B)\nearrow f.$ By Theorem \ref{t1.2}, there exist $\widetilde{u_j},\widetilde{v_j}\in PSH\cap C(\overline{B})$ such that 
$$det(\widetilde{u_j})=f_jdV
 \,\,\,~~~\text{and}~~~det(\widetilde{v_j})=g_jdV.$$ 
By the comparison principle we have $\widetilde{u_j}\searrow u,\widetilde{v_j}\searrow v$. Then (\ref{1.4}) also holds by Lemma \ref{l3.2}.\endproof

\section{Subsolution theorem of quaternionic Monge-Amp\`{e}re equation}
The purpose of this section is to prove Theorem \ref{t1.5}. Here we used the method from Ko\l odziej \cite{kolodziej94,kolodziej95,kolodziej05} for the complex Monge-Amp\`{e}re equation. Nguyen \cite{subsolution} also use Ko\l odziej's method to study the Dirichlet problem for the complex Hessian equation.\vskip 3mm

\emph{Proof of Theorem \ref{t1.5}} First we can assume that $\mu$ has compact support in $\Omega$. This is because, for non-compactly supported measure $\mu$, we can take a nondecreasing sequence of cut-off functions $\chi_j\nearrow 1$ on $\Omega$. Then $\chi_j\mu$ have compact support in $\Omega$. By Lemma \ref{comparison}, the solutions corresponding to $\chi_j\mu$ will be bounded from below by the given subsolution and they will decrease to the solution for $\mu$ by Lemma \ref{l3.2}.

Next we can modify the subsolution $v$ such that $v$ is PSH in a neighborhood of $\overline{\Omega}$ and $\lim_{\zeta\rightarrow q}v(\zeta)=0$ for any $q\in\partial\Omega$. Take an open subset $U$ such that supp$\mu\Subset U\Subset \Omega$ and define the PSH envelope
$$\widetilde{v}=\sup\{w\in PSH(\Omega):w\leq0,w\leq v~\text{in}~ U\}.$$ By Proposition \ref{p1.1}, $\widetilde{v}=\widetilde{v}^*\in PSH(\Omega)$. And $v=\widetilde{v}$ in $U$, $\lim_{\zeta\rightarrow q}\widetilde{v}(\zeta)=0$ for any $q\in\partial\Omega$. Since supp$\mu\Subset U$, we have $(\triangle \widetilde{v})^n\geq d\mu$. Take a defining function $\rho$ of $\Omega$ which is smooth and strictly PSH in a neighborhood $\Omega_1$ of $\Omega$. Since $\widetilde{v}$ is bounded, we can assume $\rho\leq \widetilde{v}$ in $\overline{U}$. Define $$\hat{v}=\begin{cases}
 \max\{\rho,\widetilde{v}\}~~~~\text{on}~~~~~~\overline{\Omega},\\
\rho~~~~~~~~~~~~~~~~\text{on}~~~\Omega_1\backslash\overline{\Omega}.
 \end{cases}$$ By Proposition \ref{p1.1}, $\hat{v}$ is a PSH function on a neighborhood of $\overline{\Omega}$ satisfying $(\triangle \hat{v})^n\geq d\mu$, $\lim_{\zeta\rightarrow q}\hat{v}(\zeta)=0$ for any $q\in\partial\Omega$. We still write $v$ instead of $\hat{v}$. Futhermore, we can make the support of $d\nu:=(\triangle v)^n$ compact in $\Omega$.
 
 Now for such subsolution $v$ we take the standard smooth regularization $w_j\searrow v$ in a neighborhood of $\overline{\Omega}$ (cf. Proposition \ref{p1.1} (5)). Let $(\triangle w_j)^n=g_jdV$. By the Radon-Nikodym theorem, $d\mu=hd\nu$ with $0\leq h\leq 1$. Denote $\mu_j=hg_jdV$. Since $w_j\searrow v$, $h(\triangle w_j)^n$ converges weakly to $h(\triangle v)^n$ by Lemma \ref{l3.2}, i.e., $\mu_j\rightarrow \mu$ as $j\rightarrow \infty$. As $\mu$ has compact support, so does $\mu_j$'s.
 Note that $hg_j\in L^p(\Omega)$ for every $p>0$. Then by Theorem \ref{t1.2}, we can find $u_j$ satisfying 
\ben\begin{cases}
 \,\,\,u_j \in PSH(\Omega)\cap C(\overline{\Omega})\\
 \,\,\,(\triangle u_j)^n= \mu_j=hg_jdV\\
 \,\,\,u_j=\varphi\,\,\,\,\text{on}\,\,\partial\Omega.
 \end{cases}\een
As we shall see the function $u:=(\limsup u_j)^*$ solves the equation (\ref{d4}). By passing to a subsequence we assume that $u_j$ converges to $u$ in $L^1(\Omega)$. Since the smooth regularization $\{w_j\}$ is uniformly bounded, we can choose a uniform $C$ such that $w_j-C<\varphi$ on $\partial\Omega$. By the comparison principle, $w_j-C\leq u_j\leq \sup_{\overline{\Omega}}\varphi$. It follows that $\{u_j\}$ is uniformly bounded, thus $u$ defined above is bounded. Now we need the following lemmas.

\begin{lem}\label{t2.4}(quasicontinuity theorem, Theorem 1.1 in \cite{wan4}) Let $\Omega$ be an open subset of $\mathbb{H}^n$ and let $u$ be a locally bounded PSH function. Then for each $\varepsilon>0$, there exists an open subset $\omega$ of $\Omega$ such that $cap(\omega)<\varepsilon$ and $u$ is continuous on $\Omega\backslash\omega$. \end{lem}
Here the quaternionic capacity of $E$ in $\Omega$ is defined in \cite{wan4} as
\begin{equation}\label{capacity}cap(E)=cap(E,\Omega)=\sup\left\{\int_{E}(\triangle
u)^n:u\in
PSH(\Omega),0<u<1\right\}.
\end{equation}

\b{lem}\label{l1} If for any $a>0$ and any compact $K\subset \Omega$ we have \beq\label{1.5}\lim_{j\rightarrow\infty}\int_{K\cap E_j(a)}(\triangle u_j)^n=\lim_{j\rightarrow\infty}\mu_j(K\cap E_j(a))=0\eeq with  $E_j(a):=\{u-u_j\geq a\}$, then the function $u$ defined above solves the Dirichlet problem (\ref{d4}).
\e{lem}
\proof By Demailly's inequality (Proposition 3.5 in \cite{wan8}) we have 
\beq\b{aligned} \label{1.6} \mu_j=(\triangle u_j)^n=&\chi_{E_j(a)}(\triangle u_j)^n+\chi_{\{u-u_j<a\}}(\triangle u_j)^n\\\leq& \chi_{E_j(a)}\mu_j+(\triangle \max\{u-a,u_j\})^n.\e{aligned}
\eeq By (\ref{1.5}), for any $s$ we can choose $j(s)$ such that $$\mu_j(K\cap E_j(\frac{1}{s}))<\frac{1}{s}, ~~j\geq j(s).$$ Let $\rho_s:=\max\{u-\frac{1}{s},u_{j(s)}\},$ then $\rho_s$ is PSH by Proposition \ref{p1.1}. By (\ref{1.6}) we have $\mu\leq \liminf_{s\rightarrow \infty}(\triangle \rho_s)^n$. By the Hartogs lemma (Theorem 4.1.9 in \cite{hormander1}), $\rho_s$ converges uniformly to $u$ on any compact $E$ such that $u|_E$ is continuous. Therefore by the quasicontinuity of quaternionic PSH functions (Lemma \ref{t2.4}) we have $\rho_s\rightarrow u$ in capacity. Thus $(\triangle \rho_s)^n\rightarrow (\triangle u)^n$ weakly by the convergence theorem (cf. Theorem 4.1 in \cite{wan7}). It follows that $\mu\leq (\triangle u)^n$. 

To show the reverse inequality, we shall prove that for $\epsilon>0$, $\mu(\Omega)\geq \int_{\Omega_\epsilon}(\triangle u)^n$, where $\Omega_\epsilon=\{q\in \Omega:\text{dist}(q,\partial\Omega)>\epsilon\}$. Note that $\rho_s=u_{j(s)}$ on a neighborhood of $\partial\Omega_\epsilon$ for $\epsilon$ small enough. By (\ref{2.37}), (\ref{2.11}) and (\ref{integral}), we can use the Stokes-type formula (Lemma \ref{p3.9}) to get
\ben\b{aligned}  &\mu(\Omega)\geq \mu(\overline{\Omega}_\epsilon)\geq \liminf_{j(s)\rightarrow\infty}\mu_{j(s)}(\Omega_\epsilon)=\liminf_{j(s)\rightarrow\infty}\int_{\Omega_\epsilon}(\triangle u_{j(s)})^n\\=&\liminf_{j(s)\rightarrow\infty}\int_{\Omega_\epsilon}d_0\left[d_1u_{j(s)}\wedge (\triangle u_{j(s)})^{n-1}\right]\\=&
\liminf_{j(s)\rightarrow\infty}\int_{\Omega_\epsilon}\sum_{i_1}\nabla_{i_10 }\left[\sum_{j_1,i_2,\ldots}
\nabla_{j_11 }u_{j(s)}\triangle_{i_2j_2}u_{j(s)}\ldots\triangle_{i_nj_n}u_{j(s)}\delta^{i_1j_1\ldots
i_nj_n}_{01\ldots(2n-1)}\right]dV
\\=&
\liminf_{j(s)\rightarrow\infty}\int_{\partial\Omega_\epsilon}
\sum_{i_1}\left[\sum_{j_1,i_2,\ldots}
\nabla_{j_11 }u_{j(s)}\triangle_{i_2j_2}u_{j(s)}\ldots\triangle_{i_nj_n}u_{j(s)}\delta^{i_1j_1\ldots
i_nj_n}_{01\ldots(2n-1)}\right]n_{i_10}dS\\=&
\liminf_{j(s)\rightarrow\infty}\int_{\partial\Omega_\epsilon}
\sum_{i_1}\left[\sum_{j_1,i_2,\ldots}
\nabla_{j_11 } \rho_s\triangle_{i_2j_2} \rho_s\ldots\triangle_{i_nj_n} \rho_s\delta^{i_1j_1\ldots
i_nj_n}_{01\ldots(2n-1)}\right]n_{i_10}dS
\\=&\liminf_{j(s)\rightarrow\infty}\int_{\Omega_\epsilon}(\triangle \rho_s)^n\geq \int_{\Omega_\epsilon}(\triangle u)^n,\e{aligned}
\een 
where $ n_{i\alpha}$ is defined by (\ref{sss4.3}) and $dS$ denotes the surface measure of $\partial\Omega$. Let $\epsilon\rightarrow 0$, we get $\mu(\Omega)\geq(\triangle u)^n(\Omega).$
\endproof

\begin{lem}\label{cln}(Chern-Levine-Nirenberg type estimate, Proposition 3.10 in \cite{wan3}) Let
$\Omega$ be a domain in $\mathbb{H}^n$. Let $K,L$ be compact subsets
of $\Omega$ such that $L$ is contained in the interior of $K$. Then
there exists a constant $C$ depending only on $K,L,\Omega$ such that for
any $v\in PSH(\Omega)$ and $u_1,\ldots u_n\in PSH\cap C^2(\Omega)$, one
has
\addtocounter{equation}{1}
 \begin{align}
 &\|\triangle
u_1\wedge\ldots\wedge\triangle u_k\|_L \leq
C\|u_1\|_{L^{\infty}(K)}\ldots\|u_k\|_{L^{\infty}(K)},
\tag{a}\label{a}\\
&\|\triangle
u_1\wedge\ldots\wedge\triangle u_k\|_L \leq
C\|u_1\|_{L^{1}(K)}\|u_2\|_{L^{\infty}(K)}\ldots\|u_k\|_{L^{\infty}(K)},\tag{$b$}\label{b}\\
&\|v\triangle
u_1\wedge\ldots\wedge\triangle u_k\|_L \leq
C\|v\|_{L^{1}(K)}\|u_1\|_{L^{\infty}(K)}\ldots\|u_k\|_{L^{\infty}(K)}.\tag{$c$}\label{c}
 \end{align}\end{lem}

\b{lem}\label{l2} Suppose that there is a subsequence of $\{u_j\}$ (we still write $\{u_j\}$) such that 
\beq\label{1.7}\int_{E_j(a_0)}(\triangle u_j)^n>A_0,~~A_0>0,a_0>0.\eeq
Then there exist $a_p>0,A_p>0,k_1>0$ such that
\beq\label{1.8}\int_{E_j(a_p)}(\triangle v_j)^{n-p}\wedge (\triangle v_k)^{p}>A_p,~~j>k>k_1,\eeq
for $v_j$'s the solutions (by Lemma \ref{l1.4}) of the Dirichlet problem \ben\begin{cases}
 \,\,\,v_j \in PSH(\Omega)\cap C(\overline{\Omega})\\
 \,\,\,(\triangle v_j)^n= g_jdV\\
 \,\,\,v_j=0\,\,\,\,\text{on}\,\,\partial\Omega.
 \end{cases}\een
\e{lem}
\proof We will prove it by induction over $p$. For $p=0$, the result (\ref{1.8}) follows from the hypothsis (\ref{1.7}) and comparison principle. Suppose that (\ref{1.8}) is true for $p<n$ and now we shall prove it for $p+1$. Note that $\{v_j\}$ is also uniformly bounded. We can assume that $-1<u_j,v_j<0$. By Chern-Levine-Nirenberg estimate (Lemma \ref{cln}) there exists $C>0$ such that \beq \label{1.9}\int_\Omega (\triangle v_j)^q\wedge (\triangle u_j)^{n-q} \leq C,\eeq for any $q=0,\ldots,n$.
 Set $$S:=(\triangle v_j)^{n-p-1}\wedge (\triangle v_k)^p.$$ By the induction hypothesis, there exist $a_p,A_p>0$ and $k_1>0$ such that 
\beq\label{2.3}\int_{E_j(a_p)}\triangle v_j\wedge S>A_p,~~j>k>k_1.\eeq
By Lemma \ref{t2.4}, for fixed $\epsilon\in (0,\frac{a_pA_p}{4(1+C)})$ (where $C$ is from (\ref{1.9})), we can choose an open set $U$ such that $cap(U,\Omega)<\frac{\epsilon}{2^{n+1}}$, and $u,v$ are continuous off the set $U$. Then 
\beq\label{2.0}\int_U(\triangle (v_j+v_k))^n<2^n cap(U,\Omega)<\frac{\epsilon}{2},~~\text{and}~~\int_U(\triangle (u_j+v_k))^n<\frac{\epsilon}{2}.\eeq
Set \ben J_1:=\int_\Omega (u-u_j)\triangle v_k\wedge S,\quad J_2:=\int_\Omega (u-u_j)\triangle v_j\wedge S.\een
Since $u=u_j=\varphi$ and $v_j=0$ on the boundary $\partial\Omega$, we can use the integration by parts formula (\ref{3.26}) to get
\beq\label{2.1}\b{aligned} J_2-J_1=&\int_\Omega v_j\triangle(u-u_j) \wedge S-\int_\Omega v_k\triangle (u-u_j)\wedge S\\=&\int_{\Omega\backslash U} (v_j-v_k)\triangle(u-u_j) \wedge S+\int_{ U} (v_j-v_k)\triangle(u-u_j) \wedge S \\ \leq&\int_{\Omega\backslash U} \|v_j-v_k\|\triangle(u+u_j) \wedge S+\int_{ U} \triangle(u+u_j) \wedge S, \e{aligned}\eeq where the last inequality follows from $-1<u_j,v_j<0$. Note that $v_j$ converges uniformly to $v$ on $\Omega\backslash U$. Let $l>k_1$ such that $ \|v_j-v_k\|<\frac{\epsilon}{4C}$ on $ \Omega\backslash U$ for $j>k>l>k_1$, where $C$ is the constant in (\ref{1.9}). By (\ref{1.9}) and (\ref{2.0}), each of the integral on the right hand side of (\ref{2.1}) does not exceed $\frac{\epsilon}{2}$. Therefore $$J_2-J_1\leq \epsilon\leq \frac{a_pA_p}{4},$$ for $j>k>l>k_1$. 

From the upper bound of all $u_j$ (resp. $v_j$) by $\sup\varphi$ (resp. $0$) on the boundary, we have for $k>k_2>l$, in a neighborhood of $\partial\Omega$ 
\beq\label{2.2} v_k\leq v+\epsilon,~~\text{and}~~u_k\leq u+\epsilon.\eeq
And they still hold on $\Omega\backslash U$ by the Hartogs lemma. By using (\ref{1.9}) (\ref{2.3}) (\ref{2.0}) and (\ref{2.2}) we get for $j>k>k_2$,
\ben\b{aligned} J_2&\geq a_p\int_{\{u-u_j\geq a_p\}}\triangle v_j\wedge S+\int_{\{u-u_j< a_p\}}(u-u_j)\triangle v_j\wedge S\\&=a_p\int_{\{u-u_j\geq a_p\}}\triangle v_j\wedge S+\int_{\{u-u_j< a_p\}\cap (\Omega\backslash U)}(u-u_j)\triangle v_j\wedge S\\&\quad+\int_{\{u-u_j< a_p\}\cap U}(u-u_j)\triangle v_j\wedge S\\ &\geq a_p\int_{\{u-u_j\geq a_p\}}\triangle v_j\wedge S-\epsilon\int_{ \Omega\backslash U}\triangle v_j\wedge S-\int_{ U}\triangle v_j\wedge S\\ &\geq a_pA_p-\epsilon(C+1)\geq \frac{3a_pA_p}{4}.
\e{aligned}\een
Fix $d>0$. It follows from (\ref{1.9}) that
$$J_1\leq \int_{\{u-u_j\geq d\}}\triangle v_k\wedge S+d\int_\Omega\triangle v_k\wedge S\leq \int_{\{u-u_j\geq d\}}\triangle v_k\wedge S+dC.$$
If we take $$a_{p+1}=d:=\frac{a_pA_p}{4C},$$ then 
\ben\b{aligned}  \int_{E_j(a_{p+1})}\triangle v_k\wedge S\geq J_1-dC\geq J_2-\epsilon-dC\geq \frac{a_pA_p}{4}:=A_{p+1},
\e{aligned}\een for $j>k>k_2$, which concludes the proof of the inductive step of Lemma \ref{l2}.\endproof

Now we continue to prove Theorem \ref{t1.5}. It suffices to show (\ref{1.5}) in Lemma \ref{l1}. Suppose that it is not true, then by the assumption of Lemma \ref{l2} we have for $p=n$ and fixed $k>k_1$, \ben\int_{E_j(a_n)} (\triangle v_k)^{n}>A_n,~~j>k.\een
Since $ (\triangle v_k)^{n}\leq M_kdV$ for some $M_k>0$, we have 
$$ V(E_j(a_n))\geq M_k^{-1}\int_{E_j(a_n)} (\triangle v_k)^{n}>\frac{A_n}{M_k},~~j>k,$$
which contradicts the fact that $u_j\rightarrow u$ in $L_{loc}^1$. Thus Theorem \ref{t1.5} follows. \qed

It follows from the subsolution theorem that if the Dirichlet problem (\ref{d4}) is solvable for a measure $\nu$ and $\mu\leq \nu$, then the Dirichlet problem (\ref{d4}) is solvable for $\mu$ as well. 
\b{cor}Let $\Omega$ be a strictly pseudoconvex domain in $\mathbb{H}^n$ and $v_1,\ldots,v_n\in PSH\cap L^\infty(\Omega)$. There exists $u\in PSH\cap L^\infty(\Omega)$ satisfying any given continuous boundary data and $(\triangle u)^n=\triangle v_1\wedge \ldots \wedge \triangle v_n$.
\e{cor}
\proof One can choose $v_1+\ldots+v_n$ as a subsolution after adding a suitable PSH function.\endproof

\b{appendix}
\section{Proof of Lemma 3.7}
\proof For $C^2$ smooth plurisubharmonic functions $u_i$, $i=1,\ldots,n$, their quaternionic Hessian $\left[\frac{\partial^2u_i}{\partial \overline{q}_j\partial q_k}(q)\right]$ are positive definite hyperhermitian matrices. By Theorem 1.1.15 in \cite{alesker1}, $$det(u_1,\ldots,u_n):=det\left((\frac{\partial^2u_1}{\partial \bar{q}_j\partial q_k}),\ldots,(\frac{\partial^2u_n}{\partial \bar{q}_j\partial q_k})\right)>0.$$
It follows from Aleksandrov inequality (cf. Corollary 1.1.16 in \cite{alesker1}, or Corollary 2.16 in \cite{zhu}) that 
\beq\label{5.1} det(u_1,u_2,\ldots,u_n)\geq det(u_1,u_1,u_3,\ldots,u_n)^{\frac{1}{2}}\cdot det(u_2,u_2,u_3,\ldots,u_n)^{\frac{1}{2}}.
\eeq
We claim that when all functions are in $PSH\cap C^2(\Omega)$, 
\begin{equation}\label{5.2}\begin{aligned}&det(\underbrace{u_1,\ldots,u_1}_p,\underbrace{u_2,\ldots,u_2}_q,v_1,\ldots,v_{n-p-q})\\\geq&det(\underbrace{u_1,\ldots,u_1}_{p+q},v_1,\ldots,v_{n-p-q})^{\frac{p}{p+q}}\cdot det(\underbrace{u_2,\ldots,u_2}_{p+q},v_1,\ldots,v_{n-p-q})^{\frac{q}{p+q}}.
\end{aligned}\end{equation}
 We prove this claim by induction. The case for $p=q=1$ holds by (\ref{5.1}). Assume by induction that the case for $p+q\leq m$ has already been proved. It suffices to prove it for $p+q\leq m+1$. First we need the following inequality.
\begin{equation}\label{5.3}det(u_1,\underbrace{u_2,\ldots,u_2}_{p+q},v,\ldots)\geq det(\underbrace{u_1,\ldots,u_1}_{p+q+1},v,\ldots)^{\frac{1}{p+q+1}}\cdot det(\underbrace{u_2,\ldots,u_2}_{p+q+1},v,\ldots)^{\frac{p+q}{p+q+1}}.
\end{equation}
By induction assumption we have
\begin{equation*}\begin{aligned}&det(u_1,\underbrace{u_2,\ldots,u_2}_{p+q},v,\ldots)\\\geq&det(\underbrace{u_1,\ldots,u_1}_{p+q},u_2,v,\ldots)^{\frac{1}{p+q}}\cdot det(\underbrace{u_2,\ldots,u_2}_{p+q},u_2,v,\ldots)^{\frac{p+q-1}{p+q}}\\=&det(\underbrace{u_1,\ldots,u_1}_{p+q-1},u_2,u_1,v,\ldots)^{\frac{1}{p+q}}\cdot det(\underbrace{u_2,\ldots,u_2}_{p+q+1},v,\ldots)^{\frac{p+q-1}{p+q}}
\\\geq&\left[det(\underbrace{u_1,\ldots,u_1}_{p+q},u_1,v,\ldots)^{\frac{p+q-1}{p+q}}\cdot det(\underbrace{u_2,\ldots,u_2}_{p+q},u_1,v,\ldots)^{\frac{1}{p+q}}\right]^{\frac{1}{p+q}}\\&\qquad\cdot det(\underbrace{u_2,\ldots,u_2}_{p+q+1},v,\ldots)^{\frac{p+q-1}{p+q}}.
\end{aligned}\end{equation*}
Then we have
\begin{equation*}\begin{aligned}&det(u_1,\underbrace{u_2,\ldots,u_2}_{p+q},v,\ldots)^{\frac{(p+q)^2-1}{(p+q)^2}}\\&\geq det(\underbrace{u_1,\ldots,u_1}_{p+q+1},v,\ldots)^{\frac{p+q-1}{(p+q)^2}}\cdot det(\underbrace{u_2,\ldots,u_2}_{p+q+1},v,\ldots)^{\frac{p+q-1}{p+q}}.
\end{aligned}\end{equation*}It follows (\ref{5.3}). Now we complete the induction by using (\ref{5.3}).
\begin{equation*}\begin{aligned}&det(\underbrace{u_1,\ldots,u_1}_{p+1},\underbrace{u_2,\ldots,u_2}_{q},v,\ldots)\\&\geq det(\underbrace{u_1,\ldots,u_1}_{p+q},u_1,v,\ldots)^{\frac{p}{p+q}}\cdot det(\underbrace{u_2,\ldots,u_2}_{p+q},u_1,v,\ldots)^{\frac{q}{p+q}}\\&\geq det(\underbrace{u_1,\ldots,u_1}_{p+q+1},v,\ldots)^{\frac{p}{p+q}}\cdot\\ &\qquad \left[det(\underbrace{u_2,\ldots,u_2}_{p+q+1},v,\ldots)^{\frac{p+q}{p+q+1}}\cdot det(\underbrace{u_1,\ldots,u_1}_{p+q+1},v,\ldots)^{\frac{1}{p+q+1}}\right]^{\frac{q}{p+q}}\\&=det(\underbrace{u_1,\ldots,u_1}_{p+q+1},v,\ldots)^{\frac{p+1}{p+q+1}}\cdot det(\underbrace{u_2,\ldots,u_2}_{p+q+1},v,\ldots)^{\frac{q}{p+q+1}}.
\end{aligned}\end{equation*} Then (\ref{5.2}) is proved.
It follows that \begin{equation*}det(u_1,\underbrace{u_2,\ldots,u_2}_{n-1})\geq (det u_1)^{\frac{1}{n}}(det u_2)^{\frac{n-1}{n}}.\end{equation*} This is the case of (\ref{5.4}) for $u_2=\ldots=u_n=u$. Assume that (\ref{5.4}) is proved for $u_{p+1}=\ldots=u_n=u$. We now prove it for $u_{p+2}=\ldots=u_n=u$. By (\ref{5.2}) we have
\begin{equation*}\begin{aligned}&det(u_1,\ldots,u_p,u_{p+1},\underbrace{u,\ldots,u}_{n-p-1})\\\geq&det(u_1,\ldots,u_p,\underbrace{u_{p+1},\ldots,u_{p+1}}_{n-p})^{\frac{1}{n-p}} \cdot det(u_1,\ldots,u_p,\underbrace{u,\ldots,u}_{n-p})^{\frac{n-p-1}{n-p}}\\\geq&\left[(det u_1)^{\frac{1}{n}}\cdots(det u_p)^{\frac{1}{n}}(det u_{p+1})^{\frac{n-p}{n}}\right]^{\frac{1}{n-p}}\cdot\left[(det u_1)^{\frac{1}{n}}\cdots(det u_p)^{\frac{1}{n}}(det u)^{\frac{n-p}{n}}\right]^{\frac{n-p-1}{n-p}}\\=&(det u_1)^{\frac{1}{n}}\cdots(det u_p)^{\frac{1}{n}}(det u_{p+1})^{\frac{1}{n}}(det u)^{\frac{n-p-1}{n}}.
\end{aligned}\end{equation*}The induction is complete. \endproof
\e{appendix}

\begin{Acknw}
This work is supported by National Nature Science Foundation in China (No.
11571305).
\end{Acknw}
\section*{References}
\bibliographystyle{plain}

 \end{document}